\documentclass[amssymb,twoside,10pt,amscd,reqno]{amsart}
\usepackage{amsmath}
\usepackage{amscd}
\usepackage{latexsym}
\usepackage{amsthm}
\usepackage{amssymb}
\theoremstyle{plain}
\newtheorem{theorem}{Theorem}[section]
\newtheorem{thm}[theorem]{Theorem}
\newtheorem{cor}[theorem]{Corollary}
\newtheorem{lem}[theorem]{Lemma}
\newtheorem{prop}[theorem]{Proposition}
\newtheorem{defn}[theorem]{Definition}

\theoremstyle{definition}

\theoremstyle{remark}

\newcommand{\mc}{\mathcal}

\newcommand{\ZZ}{\mathbb{Z}}

\newcommand{\PP}{\mathbb{P}}

\newcommand{\SP}{\text{Spec }}

\def\cG{{\mc G}}
\def\cX{{\mc X}}
\def\cE{{\mc E}}
\def\cF{{\mc F}}

\def\cK{{\mc K}}

\newcommand{\lt}{\left}
\newcommand{\rt}{\right}

\newcommand{\OO}{\mathcal O}

\begin{document}

\title{Quot functors for Deligne-Mumford stacks}
\author[Olsson]{Martin Olsson}

\author[Starr]{Jason Starr}

\maketitle

\noindent \small
\address{Department of Mathematics,  Massachusetts Institute of
Technology,  Cambridge, MA 02139,} \email{molsson@math.mit.edu}

\

\noindent
\address{Department of Mathematics, Massachusetts Institute of
Technology, Cambridge, MA 02139,} \email{jstarr@math.mit.edu}

\

\normalsize  \centerline {\emph{Dedicated to S. Kleiman on the
occasion of his 60th birthday}}

\begin{abstract}
  Given a separated and locally finitely-presented Deligne-Mumford
  stack $\cX$ over an algebraic space $S$, and a locally
  finitely-presented $\OO_{\cX}$-module $\cF$, we prove that the Quot
  functor $\text{Quot}(\cF/\cX/S)$ is represented by a separated and
  locally finitely-presented algebraic space over $S$.  Under
  additional hypotheses, we prove that the connected components of
  $\text{Quot}(\cF/\cX/S)$ are quasi-projective over $S$.
\end{abstract}


\tableofcontents

\section{Statement of results}

Let $p:\cX\rightarrow S$ be a separated, locally finitely-presented
$1$-morphism from a Deligne-Mumford stack $\cX$ to an algebraic space
$S$.  Let $\cF$ be a quasi-coherent $\OO_\cX$-module (on the \'etale
site of $\cX$) such that $\cF$ is locally finitely-presented.  Define
a contravariant functor
\begin{equation}
Q=Q(\cF/\cX/S): S-\text{schemes} \rightarrow \text{Sets}
\end{equation}
as follows.  For each $S$-scheme $f:Z\rightarrow S$, define $\cX_Z$ to
be $\cX\times_S Z$, and define $\cF_Z$ to be the pullback of $\cF$ to
$\cX_Z$.  We define $Q(Z)$ to be the set of $\OO_{\cX_Z}$-module
quotients $\cF_Z\rightarrow G$ which satisfy
\begin{enumerate}
\item $G$ is a quasi-coherent $\OO_{\cX_Z}$-module which is locally
  finitely-presented,
\item $G$ is flat over $Z$,
\item the support of $G$ is proper over $Z$.
\end{enumerate}
All of these properties are preserved by base-change on $Z$, and
therefore pullback makes $Q$ into a contravariant functor.

\begin{thm}~\label{thm-thm1}  $Q$ is represented by an algebraic space
  which is separated and locally finitely-presented over $S$.  If
  $\cF$ has proper support over $S$, then $Q\rightarrow S$ satisfies
  the valuative criterion for properness.
\end{thm}

\textbf{Remark:} Under the hypothesis that $\cF$ has proper support
over $S$, we are not claiming that $Q\rightarrow S$ is proper, because
we do not show that $Q\rightarrow S$ is quasi-compact (in general it
need not be).

\

We give a better description of $Q$ under additional hypotheses on
$\cX$.  Our first hypothesis is that $\cX$ is a global quotient.

\begin{defn}~\label{defn-gq}  Let $S$ be an algebraic space.  A
  \emph{global quotient stack} over $S$ is an (Artin) algebraic stack
  $\cX$ over $S$ which is isomorphic to a stack of the form $[Z/G],$
  where $Z$ is an algebraic space which is finitely-presented over
  $S$, and $G$ is a flat, finitely-presented group scheme over $S$
  which is a subgroup scheme of the general linear group scheme
  $\text{GL}_{n,S}$ for some $n$.
\end{defn}

This is essentially ~\cite[definition 2.9]{EHKV} (with the Noetherian
hypotheses replaced by a finite-presentation hypothesis). We remind
the reader of the following characterization in~\cite[remark
2.11]{EHKV}.

\begin{lem}~\label{rmk-gq}
  Suppose that $\cX$ is a global quotient stack over $S$ which is
  isomorphic to $[Z/G]$ as above.  Then the diagonal action of $G$ on
  $Z\times_S GL_{n,S}$ is free and the quotient $Z'=[Z\times_S
  GL_{n,S}/G]$ is an algebraic space with a right action of
  $GL_{n,S}$.  The quotient stack $[Z'/GL_{n,S}]$ is isomorphic to the
  original stack $[Z/G]$.
\end{lem}

So every quotient stack is isomorphic to a stack of the form
$[Z'/GL_{n,S}]$.

\

Our second hypothesis is that $\cX$ is a \emph{tame} Deligne-Mumford
stack.

\begin{defn}~\label{defn-tame}
  A Deligne-Mumford stack $\cX$ is \emph{tame} if for any
  algebraically closed field $k$ and any $1$-morphism $\zeta:\SP
  k\rightarrow \cX$, the stabilizer group $\text{Aut}(\zeta)(\SP k)$
  has order prime to char$(k)$, where $\text{Aut}(\zeta)$ is the
  finite $k$-group scheme defined to be the Cartesian product of the
  diagram
\begin{equation}
\begin{CD}
  @. \SP k\\
  @. @VV(\zeta, \zeta )V \\
  \cX @>\Delta >> \cX \times \cX.
\end{CD}
\end{equation}
Here $\Delta:\cX\rightarrow \cX \times \cX$ is the diagonal morphism.
\end{defn}

\begin{thm}~\label{thm-thm2}
  Suppose that $S$ is an affine scheme, and let $f:\cX\rightarrow S$
  be a separated $1$-morphism from a tame Deligne-Mumford stack to $S$
  such that $\cX$ is a global quotient over $S$ and such that the
  coarse moduli space $X$ of $\cX$ is a quasi-projective $S$-scheme
  (resp. projective $S$-scheme).  Then the connected components of $Q$
  are quasi-projective $S$-schemes (resp.  projective $S$-schemes).
\end{thm}

\textbf{Remark:} The existence of the coarse moduli space for $\cX$
follows from ~\cite{K-M}.

\

\textbf{Remark:} The condition that $S$ be an affine scheme is
required because the property of being quasi-projective is not Zariski
local on the base.

\

\textbf{Acknowledgments:}  The authors benefitted from useful
conversations with Johan de Jong.  The authors also
wish to thank the referee,
particularly for his suggestion to simplify section~\ref{gensheavsec}
by including theorem~\ref{thm-genpt}.

This work was largely done
during the semester program on Algebraic Stacks, Intersection Theory
and Non-Abelian Hodge Theory at MSRI during the spring of 2002
where both authors were MSRI
Postdoctoral Fellows.  The first author was partially supported by an NSF
Postdoctoral Fellowship.
The second author was partially supported by NSF grant DMS-0201423.

\section{Representability by an algebraic space}~\label{sec-M}




In this section we prove theorem \ref{thm-thm1}.

Note first that $Q$ is a sheaf for the fppf-topology by descent
theory, and is limit preserving. In addition, for each open
substack $\mc U\subset \cX$ there is a natural open immersion
$Q(\mc F|_{\mc U}/\mc U/S)\subset Q$.  Moreover, $Q$ is the union
over finitely presented open substacks $\mc U\subset \mc X$ of the
$Q(\mc F|_{\mc U}/\mc U/S)$.  We may therefore assume that $\mc X$
is of finite presentation over $S$. Since the question of
representability of $Q$ is \'etale local on $S$ we may assume that
$S$ is an affine scheme, and by a standard limit argument we may
assume that $S$ is of finite type over $\text{Spec}(\mathbb{Z})$.

Under these assumptions we
 prove the theorem by verifying the conditions of theorem 5.3 of
\cite{Artin}.

\

\noindent \emph{Commutation with inverse limits.}
We need the Grothendieck existence theorem for Deligne-Mumford stacks:
\begin{prop}\label{stackgrothex}
Let $A$ be a complete noetherian local ring, $\cX /A$  a
Deligne-Mumford stack of finite type, and let $A_n =
A/\mathfrak{m}_A^{n+1}$, $\cX_n = \cX\otimes _AA_n$. Then the
natural functor
\[
\begin{CD}
(\hbox{coherent sheaves on $\cX$ with support proper over $A$})\\
@VVV \\
(\hbox{compatible families of coherent sheaves on the $\cX_n$ with
proper support})
\end{CD}
\]
is an equivalence of categories.
\end{prop}
\begin{proof}
This proposition is perhaps best viewed in a context of formal
algebraic stacks, but since we do not want to develop such a
theory here we take a more ad-hoc approach.

Let $\widehat {\cX}$ be the ringed topos $(\cX_{0,et}, \mc O_{\widehat
{\cX}})$, where $\mc O_{\widehat {\cX}}$ is the sheaf of rings which
to any \'etale $\cX_0$-scheme $U_0\rightarrow \cX_0$ associates
\[
\varprojlim \Gamma (U_n, \mc O_{U_n}).
\]
Here $U_n$ denotes the unique lifting of $U_0$ to an \'etale
$\cX_n$-scheme.  There is a natural morphism of ringed topoi
\[
j:\widehat \cX\rightarrow \cX_{et}.
\]
If $\mc F$ is a sheaf of $\mc O_\cX$-modules, then we denote by
$\widehat {\mc F}$ the sheaf $j^*\mc F$.  Note that the functor $\mc F\mapsto
\widehat {\mc F}$ is an exact functor.

\begin{lem}\label{cohlem}
If $\mc F$ and $\mc G$ are coherent sheaves on $\cX$ with proper
support over $A$, then for every integer $n$ the natural map
\[
\hbox{Ext}^n_{\mc O_\cX}(\mc F, \mc G)\longrightarrow {Ext}^n_{\mc
O_{\widehat {\cX}}}(\widehat {\mc F}, \widehat {\mc G})
\]
is an isomorphism.
\end{lem}
\begin{proof}
Observe first that the natural map
\[
\mc Ext^n_{\mc O_\cX }(\mc F, \mc G)^{\wedge }\rightarrow \mc
Ext^n_{\mc O_{\widehat \cX }}(\widehat {\mc F},\widehat  {\mc G})
\]
is an isomorphism.
Indeed this can be verified locally and so follows from the theory of formal
schemes.
From this and the local-to-global spectral sequence for $\hbox{Ext}$ it
follows that it suffices to show that for any coherent sheaf $\mc
F$ with proper support, the natural map
\begin{equation}\label{cohmap}
H^n(\cX, \mc F)\longrightarrow H^n(\widehat {\cX}, \widehat {\mc F})
\end{equation}
is an isomorphism for every $n$.  We prove this by induction on
the dimension of $\cX$.

The key observation is that if $a:U\rightarrow \cX$ is a finite
morphism from a scheme to $\cX$, then the map \ref{cohmap} is known
to be an isomorphism for all $n$ in the case when $\mc F$ is equal
to $a_*\mc F'$ for some coherent sheaf $\mc F'$ on $U$ (since $a$
is finite and \cite{EGA}. III.5).

If the dimension of $\cX$ is zero, then we can by (\cite{LM-B},
16.6) find a finite \'etale cover $a:U\rightarrow \cX$ of $\cX$ by a
scheme $U$. If $b:U\times _\cX U\rightarrow \cX$ denotes the canonical
map, then the sequences
\[
\mc F\rightarrow a_*a^*\mc F\rightrightarrows b_*b^*\mc F
\]
\[
\widehat {\mc F}\rightarrow a_*a^*\widehat {\mc
F}\rightrightarrows b_*b^*\widehat {\mc F}
\]
are exact so the result is true for $n = 0$.

To prove the result for general $n$, assume the result is true for
$n-1$ and let $\mc G = a_*a^*\mc F/\mc F$ so that we have an exact
sequence of coherent sheaves
\[
0\rightarrow \mc F \rightarrow a_*a^*\mc F\rightarrow \mc
G\rightarrow 0.
\]
Then the commutative diagram
\[
\begin{matrix}
H^{n-1}(\cX, a_*a^*\mc F)& \rightarrow & H^{n-1}(\cX, \mc G) &
\rightarrow
& H^n(\cX, \mc F) & \rightarrow &H^n(\cX, a_*a^*\mc F)\\
\downarrow && \downarrow && \downarrow && \downarrow \\
H^{n-1}(\widehat \cX, a_*a^*\widehat {\mc F})& \rightarrow &
H^{n-1}(\widehat \cX, \widehat {\mc G}) & \rightarrow & H^n(\widehat
\cX, \widehat {\mc F}) & \rightarrow &H^n(\widehat \cX, a_*a^*\widehat
{\mc F})
\end{matrix}
\]
shows that the map $H^n(\cX, \mc F)\rightarrow H^n(\widehat \cX,
\widehat {\mc F})$ is injective.  But then the map $H^n(\cX, \mc
G)\rightarrow H^n(\widehat {\cX}, \widehat {\mc G})$ is also
injective, and an analysis of the diagram
\[
\begin{matrix}
H^{n-1}(\cX, \mc G)& \rightarrow & H^{n}(\cX, \mc F) & \rightarrow
& H^n(\cX, a_*a^*\mc F) & \rightarrow &H^n(\cX, \mc G)\\
\downarrow && \downarrow && \downarrow && \downarrow \\
H^{n-1}(\widehat \cX, \widehat {\mc G})& \rightarrow &
H^{n}(\widehat \cX, \widehat {\mc F})& \rightarrow & H^{n}(\widehat
\cX, a_*a^*\widehat {\mc F}) & \rightarrow & H^n(\widehat \cX,
\widehat {\mc G})
\end{matrix}
\]
reveals that the map $H^n(\cX, \mc F)\rightarrow H^n(\widehat \cX,
\widehat {\mc F})$ is an isomorphism.  This completes the proof of
the case when $dim(\cX) = 0$.

To prove the result for general $\cX$ we assume the result is true
for $dim(\cX)-1$ and  proceed by induction on $n$. By (\cite{LM-B},
16.6) there exists a finite surjective morphism $a:U\rightarrow \cX$
which is generically \'etale.  Let $b:U\times _{\cX }U\rightarrow \cX$ be
the natural map and let $K$ be the kernel of $a_*a^*\mc
F\rightarrow b_*b^*\mc F$. Note that the map \ref{cohmap} for $K$
is an isomorphism when $n=0$. There is a natural map $\mc
F\rightarrow K$ which is generically an isomorphism.  Let $\mc F'$
be the kernel of this map and let $\mc F^{\prime \prime }$ be the
cokernel. Then $\mc F'$ and $\mc F^{\prime \prime }$ have
lower-dimensional support and so the map \ref{cohmap} is an
isomorphism for these sheaves.  If $\mc G$ denotes the image of
$\mc F\rightarrow K$, then the map $H^0(\cX, \mc G)\rightarrow
H^0(\widehat {\cX}, \widehat {\mc G})$ is an isomorphism by the
corresponding result for $K$ and $\mc F^{\prime \prime }$ and the
exact sequence
\[
0\rightarrow \mc G\rightarrow K\rightarrow \mc F^{\prime \prime
}\rightarrow 0.
\]
Then from the exact sequence
\[
0\rightarrow \mc F'\rightarrow \mc F\rightarrow \mc G\rightarrow 0
\]
we deduce that the map $H^0(\cX, \mc F)\rightarrow H^0(\widehat \cX,
\widehat {\mc F})$ is an isomorphism.

To prove the result for $n$ assuming the result for $n-1$, note
that a similar argument to the one above shows that if  the map
$H^n(\cX, K)\rightarrow H^n(\widehat {\cX}, \widehat K)$ is injective
(resp. an isomorphism) then the map \ref{cohmap} is injective
(resp. an isomorphism).  Therefore the proof is completed by using
the argument of the case $dim(\cX) = 0$ with $K$ replacing $\mc F$.
\end{proof}

Now observe that the category of compatible families of coherent sheaves on the
$\cX_n$ with proper support is naturally viewed
as a full subcategory of the category of
sheaves of $\mc O_{\widehat \cX}$-modules on $\widehat {\cX}$ using
the functor which sends a family $\{\mc F_n\}$ to the sheaf
$\widehat {\mc F}$ associated to the presheaf
\[
U_0\mapsto \varprojlim \Gamma (U_0, \mc F_n).
\]
The functor $\{\mc F_n\}\mapsto \widehat {\mc F}$ is fully
faithful and identifies the category of compatible systems of coherent
sheaves on the $\cX_n$ with proper support with a subcategory of the
category of sheaves of $\mc O_{\widehat {\cX}}$-modules which is
closed under the formation of kernels, cokernels, and extensions.
Indeed, these assertions can be verified locally on $\cX$ and hence follow from
the corresponding statements for formal schemes.  From this it and the
lemma it follows that the functor in \ref{stackgrothex} is fully faithful.

Now the
functor from coherent sheaves on $\cX$ with proper support to
sheaves of $\mc O_{\widehat \cX}$-modules on $\widehat \cX$ identifies
the category of such sheaves with a full subcategory of the
category of sheaves of $\mc O_{\widehat \cX}$-modules which is
stable under the formation of kernels, cokernels, and extensions
(by lemma \ref{cohlem}).  Thus the following two lemmas prove
proposition \ref{stackgrothex}.

\begin{lem}
Let $\mc A$ be an abelian category and let $\mc A'\subset \mc A$
be a full subcategory which is stable under the formation of
kernels, cokernels, and extensions.  Then any object of $\mc A$
which  admits a morphism to an object of $\mc A'$ such that the
kernel and cokernel are in $\mc A'$ is in $\mc A'$.
\end{lem}
\begin{proof}
By assumption such an object $\mc F\in \mc A$ sits in an exact
sequence
\[
0\rightarrow K\rightarrow \mc F\rightarrow \mc F'\rightarrow
Q\rightarrow 0
\]
where $K, \mc F', Q\in \mc A'$.  Let $K' = Ker(\mc F'\rightarrow
Q)$.  Then $K'\in \mc A'$, and we have an exact sequence
\[
0\rightarrow K\rightarrow \mc F\rightarrow K'\rightarrow 0.
\]
Therefore $\mc F\in \mc A'$.
\end{proof}

\begin{lem}
For every compatible family of coherent sheaves $\{\mc F_n\}$ with
proper support and associated sheaf $\widehat {\mc F}$ on
$\widehat {\cX}$, there exists a morphism $\widehat {\mc
F}\rightarrow j^*\mc G$ for some coherent sheaf $\mc G$ on $\cX$
with proper support such that the kernel and cokernel are
isomorphic to the pullbacks of coherent sheaves on $\cX$ with proper
support.
\end{lem}
\begin{proof}
We proceed by induction on the dimension of the support of
$\widehat {\mc F}$.

If the support of $\widehat {\mc F}$ has dimension zero, then for
any \'etale cover $a:U\rightarrow \cX$ the inverse image of the
support of $\widehat {\mc F}$ in $\widehat {U}$ is a closed
subscheme which is of dimension zero, hence is proper.  Therefore,
there exists  a unique sheaf on $U$ inducing the restriction of
$\widehat {\mc F}$ to $\widehat {U}$ (by \cite{EGA} III.5).
Moreover, by the uniqueness this sheaf comes with descent datum
relative to $a$. Hence $\widehat {\mc F}$ is induced from a
coherent sheaf on $\cX$ with proper support.

As for the general case, choose a morphism $a:U\rightarrow \cX$
which is finite and generically \'etale (such a morphism exists by
\cite{LM-B}, 16.6), and let $b:U\times _{\cX}U\rightarrow \cX$ be the
natural map.  Then $a_*a^*\widehat {\mc F}$ and $b_*b^*\widehat
{\mc F}$ are obtained from coherent sheaves on $\cX$, and hence so
is
\[
K := Ker(a_*a^*\widehat {\mc F}\rightarrow b_*b^*\widehat {\mc
F}).
\]
Moreover, there is a natural map $\widehat {\mc F}\rightarrow K$
which is generically an isomorphism.  Hence the kernels and
cokernels have lower dimensional support and by induction are
obtained from coherent sheaves on $\cX$.
\end{proof}

\end{proof}

\noindent \emph{Separation Conditions.}  These follow by the same
reasoning as in (\cite{Artin}, page 64).

\

\noindent \emph{Deformation theory.}

Suppose given a deformation situation (in the sense of
\cite{Artin})
\[
A'\rightarrow A\rightarrow A_0
\]
and a quotient $\mc F_A\rightarrow \mc G_A$ over $\cX_A =
\cX\times _S\text{Spec}(A)$ giving an element of $Q(A)$.  Suppose
further that a map $B\rightarrow A$ is given.  Then it is
well-known (see for example \cite{Sch}, 3.4) that the map
\begin{equation}\label{schmap}
Q(A'\times _BA)\longrightarrow Q(A')\times _{Q(A)}Q(B)
\end{equation}
is a bijection.

If $M = Ker(A'\rightarrow A)$, then it follows from the
bijectivity of \ref{schmap} that a deformation theory in the sense
of (\cite{Artin}) is provided by the module $Q_{\mc G_0}(A_0[M])$
(see for example loc. cit. page 47).  Here $A_0[M]$ denotes the
ring with underlying module $A_0\oplus M$ and multiplication
\[
(a, m)\cdot (a', m') = (aa', am'+a'm),
\]
and $Q_{\mc G_0}(A_0[M])$ denotes the set of elements in
$Q(A_0[M])$ whose image in $Q(A_0)$ is the reduction $\mc G_0$ of
$\mc G$.  The conditions on the deformation theory of
(\cite{Artin}, theorem 5.3) are therefore satisfied by the
following lemma and standard properties of cohomology.
\begin{lem}
Let $\mc H_0 = Ker(\mc F_0\rightarrow \mc G_0)$ (where $\mc F_0$
denotes the reduction of $\mc F$ to $\cX_{A_0}$). Then there is a
natural $A_0$-module isomorphism
\[
Q_{\mc G_0}(A_0[M]) \simeq \hbox{Ext}^0(\mc H_0, \mc G_0\otimes
M).
\]
\end{lem}
\begin{proof}
Let
\[
\begin{CD}
0@>>> \mc G_0\otimes M@>>> \mc G_0\otimes M\oplus \mc F_0@>>> \mc F_0@>>> 0
\end{CD}
\]
be the sequence obtained by pushing out the sequence
\[
\begin{CD}
0@>>> \mc F_0\otimes M@>>> \mc F_0\oplus \mc F_0\otimes M@>>> \mc G_0@>>> 0
\end{CD}
\]
via the given map $\mc F_0\rightarrow \mc G_0$.  Then to give a lifting of
$\mc G_0$ is equivalent to giving a sub-$\mc O_{\cX _{A_0}}$-module
$\mc H\subset \mc F_0\oplus \mc G_0\otimes M$ such that the induced map to $\mc F_0$
induces an isomorphism $\mc H\simeq \mc H_0$.  Note that such a sub-module
is automatically a sub-$\mc O_{\cX _{A_0[M]}}$-module.
In other words, the set of liftings
of $\mc G_0$ is in natural bijection with the set of maps $\mc H_0\rightarrow \mc F_0
\oplus \mc G_0\otimes M$ lifting the inclusion $\mc H_0\subset \mc F_0$.  But to give
such a map is precisely equivalent to giving a morphism $\mc H_0\rightarrow \mc G_0
\otimes M$.  The verification that this bijection is a module homomorphism
is left to the reader.

\end{proof}

\noindent \emph{Conditions on the obstructions.}

Our understanding of the obstruction theory of $Q$ will be in a
2-step approach (correcting a mistake in \cite{Artin}).  Let
$A'\rightarrow A\rightarrow A_0$ be a deformation situation as
above and $\mc F_A\rightarrow \mc G_A$ an object in $Q(A)$.
For any quotient $\epsilon :M\rightarrow M_\epsilon $, let $A_\epsilon $ be the quotient
of $A'$ by the kernel of $\epsilon $.
For such an $\epsilon $,
the
first obstruction to lifting $\mc G_A$ to $A_\epsilon $ is that the map
\[
M\otimes \mc F_{A_0}\longrightarrow M_\epsilon \otimes \mc G_{A_0}
\]
factors through $M\mc F_{A'}$.  If we let $\mc T$ be the kernel
of $M\otimes \mc F_{A_0}\rightarrow \mc F_{A'}$, then we want that the map
\begin{equation}\label{tormap}
\mc T\longrightarrow M_\epsilon \otimes \mc G_{A_0}
\end{equation}
is zero.  If this is the case, then there is a canonical map $M\mc
F_{A'}\rightarrow M_\epsilon \otimes \mc G_{A_0}$ and the condition that
there exists an element in $Q(A_\epsilon )$ inducing $\mc G$ is equivalent
to the statement that the resulting extension
\[
0\rightarrow M_\epsilon \otimes \mc G_{A_0}\rightarrow E\rightarrow \mc
F_A\rightarrow 0
\]
is obtained from an extension of $\mc G_{A}$ by $M_\epsilon \otimes \mc
G_{A_0}$. In other words, if the map \ref{tormap} is zero, then if
we let $\mc H =Ker(\mc F_A\rightarrow \mc G_A)$ there exists a
canonical obstruction in
\[
\hbox{Ext}^1_{\cX_{A'}}(\mc H, M_\epsilon \otimes \mc G_{A_0})
\]
whose vanishing is necessary and sufficient for the existence of a
lifting.

Now the only condition on the obstructions in (\cite{Artin}, 5.3)
which does not follow immediately from the bijectivity of
\ref{schmap}, is condition (5.3, [5${}'$].c).

Thus suppose given a deformation situation as before, $\xi \in
Q(A)$, and suppose further that $M$ is free of rank $n$.  Let $K$
be the field of fractions of $A_0$, and denote by a subscript $K$
the localizations at the generic point of $\text{Spec}(A_0)$.  We
suppose that for every one-dimensional quotient $M_K\rightarrow
M_K^*$ there is a non-trivial obstruction to lifting $\xi _K$ to
$Q(A_K^*)$, where $A_K^*$ denotes the extension defined by
$M_K^*$.  Then we have to show that there exists a non-empty open
subset $U\subset \text{Spec}(A_0)$ such that for every quotient
$\epsilon $ of $M$ of length one with support in $U$, $\xi $ does
not lift to $A_\epsilon $ (the extension obtained from $\epsilon
$).

Let
\[
\phi \in \hbox{Ext}^0(\mc T, M\otimes \mc
G_{A_0})\simeq \hbox{Ext}^0(\mc T, \mc
G_{A_0})\otimes M
\]
be the class defined by \ref{tormap} in the case when $\epsilon $ is the identity.
 We reduce to the case when
$\phi  = 0$.  Once this reduction is made the argument of
(\cite{Artin}, page 66) will finish the proof.

To make the reduction, we can by shrinking on $\text{Spec}(A_0)$
assume that
$$ \hbox{Ext}^0(\mc T, \mc G_{A_0})$$ is a free module; say of
rank $r$.  In addition, by the argument of (\cite{Artin}, page
66), we can after shrinking $\text{Spec}(A_0)$ assume that for
each point $s\in \text{Spec}(A_0)$, the natural map
\begin{equation}\label{pointeq}
\hbox{Ext}^0(\mc T, \mc G_{A_0})\otimes
k(s) \longrightarrow  \hbox{Ext}^0(\mc T,
\mc G_{A_0}\otimes k(s))
\end{equation}
is an isomorphism.

Choosing a basis for $ \hbox{Ext}^0(\mc T,
\mc G_{A_0})$ we can think of $\phi $ as an element
\[
\phi = (\phi _1, \dots, \phi _r)\in M^r.
\]
Let $N\subset M$ be the submodule of $M$ generated by the $\phi
_i$, and let $M' = M/N$.  After further shrinking
$\text{Spec}(A_0)$ we can assume that $M'$ is a free module.  Now
note that any length-one quotient $M\rightarrow M_{\epsilon }$ for
which the obstruction to lifting $\xi $ goes to zero factors
through $M'$ by \ref{pointeq}.
 Moreover, any such quotient which does not factor through $M'$ is
obstructed. Therefore, we may replace $M$ by $M'$ and hence are
reduced to the case when $\phi = 0$.

\

\noindent \emph{Valuative criteria for properness when $\mc F$ has
proper support}

Let $R$ be a discrete valuation ring with field of fractions $K$,
and let $i:\cX_K\hookrightarrow \cX$ be the inclusion of the generic
fiber.  We suppose that we have a flat quotient $\mc
F_K\rightarrow \mc G_K$ over the generic fiber which we wish to
extend to $\cX$.  For this we take $\mc G$ to be the image of the
map $\mc F\rightarrow i_*\mc G_K$.  The image is evidently a
coherent sheaf, and has proper support since $\mc F$ has  proper support.  It
is flat because it is a subsheaf of a torsion free sheaf, and by
definition $\mc G$ induces $\mc G_K$ on the generic fiber.

This completes the proof of theorem \ref{thm-thm1}. \qed

\section{Flattening stratifications}~\label{sec-FS}

As a corollary of the representability result in the last section,
when the locally finitely-presented sheaf $\mc F$ has proper support
over $S$, we can construct the \emph{flattening stratification} of
$\mc F$ as an algebraic space.  But in fact this algebraic space is
representable and quasi-affine over $S$.  This is a crucial step in
the proof of theorem~\ref{thm-thm2}.  Therefore we include a proof of
the following fact:

\begin{prop}~\label{prop-qaff}  Suppose $f:Y\rightarrow X$ is a
  finite-type, separated, quasi-finite morphism of algebraic spaces.
  Then f is quasi-affine, in particular $f$ is representable by
  schemes.
\end{prop}

\begin{proof}
  We need to prove that the natural morphism of $X$-schemes,
\begin{equation}
\iota:Y\rightarrow \SP_X f_*\OO_Y
\end{equation}
is an open immersion.  By ~\cite[proposition II.4.18]{Kn}, $f_*\OO_Y$
commutes with flat base change on $X$, so the formation of $\iota$
commutes with flat base change on $X$.  Moreover one may check that a
morphism of $X$-schemes is an open immersion after fpqc base change.
Thus, without loss of generality, we may suppose that $X$ is an affine
scheme.

\

To prove $\iota$ is an open immersion, it suffices to check that for
each point $p\in Y$, the following two conditions hold:
\begin{enumerate}
\item $\iota$ is \'etale at $p$, and
\item there is an open set $p\in U\subset Y$ which is disjoint from
the image of $Y\times_{\iota,\iota} Y - \Delta(Y)$ under the
projection $\text{pr}_1:Y\times_{\iota,\iota} Y \rightarrow Y$.
\end{enumerate}
Here $Y\times_{\iota,\iota} Y$ is the fiber product of $Y$ with itself
over $\SP_X f_*\OO_Y$ and $\Delta:Y\rightarrow Y\times_{\iota,\iota}
Y$ is the diagonal morphism.  It is clear that if both $(1)$ and $(2)$
hold for each point $p\in Y$, then $\iota$ is an \'etale monomorphism,
and therefore $\iota$ is an open immersion.

\

The claim is that for given $p\in Y$, we may check $(1)$ and $(2)$
after passing to an \'etale neighborhood of $q=f(p)\in X$, i.e. if
$X'\rightarrow X$ is an \'etale morphism, $q'\in X'$ is a point lying
over $q$, and if $p'\in Y':= X'\times_X Y$ is the point lying over
$q'$ and $p$, then it suffices to check $(1)$ and $(2)$ for $p'$.  We
have mentioned that the natural morphism of $X'$-schemes
\begin{equation}
\iota':Y'\rightarrow \SP_{X'} (f')_*\OO_{Y'}
\end{equation}
is the base-change of $\iota$.  The property of being \'etale at a
point can be checked after \'etale (and even flat) base-change, so
if $(1)$ holds for $p'$ then $(1)$ holds for $p$.  Suppose $(2)$
holds for $p'$ and let $U'$ be an open set $p'\in U'\subset Y'$ as
in $(2)$. Let $p\in U\subset Y$ be the open image of $U'$ under
$g:Y'\rightarrow Y$.  We have the equality
\begin{eqnarray}
g(U') \cap \text{pr}_1\lt(Y\times_{\iota,\iota} Y - \Delta(Y)\rt) =
g\lt( U' \cap g^{-1} \text{pr}_1 \lt( Y\times_{\iota,\iota} Y -
\Delta(Y)\rt)\rt) = \\
g\lt( U' \cap \text{pr}'_1 \lt(
Y'\times_{\iota',\iota'} Y' - \Delta'(Y')\rt) \rt).
\end{eqnarray}
But $U'\cap \text{pr}'_1 \lt( Y'\times_{\iota',\iota'} Y' -
\Delta'(Y') \rt) = \emptyset$ by assumption.  Thus $U=g(U')$ satisfies
condition $(2)$ for $p\in Y$.  So it suffices to check $(1)$ and $(2)$
for $p'\in Y'$.

\

Now let $Z\rightarrow Y$ be an \'etale morphism of a scheme $Z$ to $Y$
and $r\in Z$ a point mapping to $p$.  Then $Z\rightarrow X$ is
finite-type, separated and quasi-finite.  By ~\cite[proposition
2.3.8(a)]{Ner}, we may find an \'etale morphism $X'\rightarrow X$ with
$X'$ an affine scheme, and a point $q'\in X'$ mapping to $q=f(p)$ such
that if $r'\in Z'$ is the point lying over $q'$ and $r$, and if $Z'_0$
is the connected component of $Z'$ containing $r'$, then
$Z'_0\rightarrow X'$ is a finite morphism of affine schemes.  Let
$Y'_0\subset Y'$ be the connected component of $Y'$ containing the
image of $Z'_0$.  Then $Z'_0\rightarrow Y'_0$ is \'etale with dense
image, but it is also proper.  Therefore $Z'_0\rightarrow Y'_0$ is
finite,\'etale and surjective.  So, by Knudsen's version of
Chevalley's theorem~\cite[theorem III.4.1]{Kn}, $Y'_0$ is an affine
scheme.

\

Let $Y'_r$ denote the union of  connected components $Y'-Y'_0$.
Then we have a decomposition of $\text{Spec}_{X'}(f')_* \OO_{Y'}$
into  components, $\text{Spec}_{X'}(f'_0)_* \OO_{Y'_0}$ union
$\text{Spec}_{X'}(f'_r)_* \OO_{Y'_r}$ and a decomposition of
$\iota'$ as the ``disjoint union'' of the two morphisms:
\begin{eqnarray}
\iota'_0:Y'_0 \rightarrow \SP_{X'}(f'_0)_* \OO_{Y'_0}, \\
\iota'_r:Y'_r \rightarrow \SP_{X'}(f'_r)_* \OO_{Y'_r}.
\end{eqnarray}
And $\iota'$ is an isomorphism.  Thus $\iota'$ is \'etale at $p'$.
And defining $U'=Y'_0$, we see that $U'$ satisfies the condition $(2)$
for $p'\in Y'$.  So $(1)$ and $(2)$ hold for $p'$ and thus $(1)$ and
$(2)$ hold for $p\in Y$.  Since this holds for every $p\in Y$, we
conclude that $\iota$ is an open immersion, i.e. $f:Y\rightarrow X$ is
quasi-affine.
\end{proof}

\textbf{Remark:} We could not find precisely this statement in
~\cite{Kn}, although it follows easily from results proved there.  If
one further assumes that $f$ is finitely-presented, then this result
also follows easily from ~\cite[th\'eor\`eme 16.5]{LM-B}.

\

Now suppose $S$ is an algebraic space, $f:\cX\rightarrow S$ is a
1-morphism from a Deligne-Mumford stack to $S$ which is separated and
locally finitely-presented.  Let $\cF$ be a locally finitely-presented
$\OO_\cX$-module.  Consider the functor
\begin{equation}
\Sigma:\lt(S-\text{schemes}\rt)^{\text opp}\rightarrow \text{Sets},
\end{equation}
which to any morphism $g:T\rightarrow S$ associates $\lt\{*\rt\}$ if
the pullback of $\cF$ to $T\times_S \cX$ is flat over $T$, and which
associates $\emptyset$ otherwise.  Given a morphism of $S$-schemes,
$h:T_1\rightarrow T_2$, the morphism $\Sigma(h)$ is defined to be the
unique map $\Sigma(T_2)\rightarrow \Sigma(T_1)$.  For this to make
sense, we must check that if $\Sigma(T_2)$ is nonempty, then so is
$\Sigma(T_1)$.  But this is clear, if the pullback of $\cF$ to
$T_2\times_S \cX$ is flat over $T_2$, then the pullback of this sheaf
over $T_1\times_{T_2} \lt(T_2 \times_S \cX\rt)$ is flat over $T_1$.
And this pullback is isomorphic to the pullback of $\cF$ to
$T_1\times_S \cX$.  So $\Sigma(T_1)$ is nonempty.

\

\begin{thm}~\label{thm-fs}
  Let $f:\cX\rightarrow S$, $\cF$, and $\Sigma$ be as above.
\begin{enumerate}
\item $\Sigma$ is an fpqc sheaf which is limit preserving and
  $\Sigma\rightarrow S$ is a monomorphism.
\item If $\cF$ has proper support over $S$, then $\Sigma$ is an
  algebraic space and $\Sigma\rightarrow S$ is a surjective,
  finitely-presented monomorphism.  In particular, $\Sigma\rightarrow
  S$ is quasi-affine.
\end{enumerate}
\end{thm}

\begin{proof}
  It is immediate that $\Sigma\rightarrow S$ is a monomorphism.  Since
  one may check that a quasi-coherent sheaf on $T\times_S \cX$ is flat
  over $T$ after performing an fpqc base change of $T$, it follows
  that $\Sigma$ is an fpqc sheaf.  The fact that $\Sigma$ is
  limit-preserving follows from ~\cite[th\'eor\`eme IV.11.2.6]{EGA}.
  So (1) is proved.

\

To prove (2), first notice that we may use (1) to reduce to the case
that $S=\SP A$ is a Noetherian affine scheme, and $\cX$ is the support
of $\cF$ which is a proper, finitely-presented Deligne-Mumford stack
over $S$.  By theorem~\ref{thm-thm1} we know the Quot functor of $\cF$
is represented by an algebraic space $Q\rightarrow S$ which is
separated and locally finitely-presented.  Denote by
\begin{equation}
\begin{CD}
  0 @>>> {\mc K} @>>> \text{pr}_2^* \cF @>>> {\mc G} @>>> 0
\end{CD}
\end{equation}
the universal short exact sequence on $Q\times_S \cX$.  Since
$\text{pr}_2^* \cF$ has proper support over $Q$,  ${\mc K}$ also
has proper support over $Q$.  Define $U$ to be the complement in $Q$
of the image of the support of ${\mc K}$.  The restriction of the
universal short exact sequence over $U$ is simply the identity
morphism of $\text{pr}_2^* \cF$ to itself, i.e. the pullback of
$\cF$ to $U$ is flat over $U$.  So we have an induced morphism
$U\rightarrow \Sigma$.  Conversely there is an obvious morphism
$\Sigma\rightarrow Q$ which factors through $U$, and we see that
$U\rightarrow \Sigma$ is a natural isomorphism of algebraic spaces.
So $\Sigma$ is an algebraic space and $\Sigma\rightarrow S$ is a locally
finitely-presented monomorphism.

\

Any module over a field is flat over that field, therefore for each
field $K$ and each morphism $\SP K \rightarrow S$, this morphism
factors through $\Sigma\rightarrow S$.  So $\Sigma\rightarrow S$ is
surjective.  The claim is that for any
surjective, locally finitely-presented monomorphism $h:Y\rightarrow S$
of an algebraic space to a Noetherian affine scheme, $Y$ is
quasi-compact.

\

We will prove this claim by induction on the dimension of $S$.
If $S=\SP K$ for some field $K$, it is obvious.  If
$Y^{\text red}$ is quasi-compact, the same is true of $Y$, so we may
reduce to the case that $S$ is reduced.  A finite union of
quasi-compact sets is quasi-compact, so we may reduce to the case that
$S$ is integral.  Now suppose the result is proved for all schemes $S$
of dimension at most $n$ and suppose $S$ is an integral scheme of
dimension $n+1$.  Let $\SP K$ be the generic point of $S$ and let
$U\subset Y$ be the open set where $U\rightarrow S$ is \'etale.  Then
$\SP K$ factors through $Y$, and in fact it is contained in $U$.  Thus
$U\rightarrow S$ is an \'etale monomorphism, i.e. an open immersion
which has dense image.  Let the complement of $U$
in $S$ be $C$ (with reduced scheme structure)
and let the preimage of $C$ be $Z$.  Then $Z\rightarrow
C$ is again a surjective, locally finitely-presented monomorphism, and
$C$ has dimension at most $n$.  So by the induction assumption, $Z$ is
quasi-compact.  Since $U$ is an open subset of a Noetherian scheme, it
is quasi-compact.  Thus $Y=U\cup Z$ is quasi-compact and the claim is
proved by induction.

\

By the last paragraph, we conclude that $\Sigma\rightarrow S$ is a
surjective, finitely-presented monomorphism of algebraic spaces.  So
by proposition~\ref{prop-qaff}, we conclude that $\Sigma \rightarrow
S$ is quasi-affine.
\end{proof}

\textbf{Remarks:} (1)  If the support of $\cF$ is a scheme and the
morphism to $S$ is projective, then it follows from ~\cite[theorem,
p.55]{Mum} that $\Sigma\rightarrow S$ is a disjoint union of locally
closed immersions.  While one can find examples of surjective,
finitely-presented monomorphisms $Y\rightarrow S$ not a disjoint union
of locally closed immersions, we conjecture that $\Sigma\rightarrow S$
is a disjoint union of locally closed immersions whenever $\cF$
has proper support over $S$.

\

(2)  Again in the case that the support of $\cF$ is projective
     over $S$, the methods in ~\cite[section 8]{Mum} provide a
     \emph{global construction} of the flattening stratification
     $\Sigma\rightarrow S$ along with a partition labelled by the
     Hilbert polynomial of the fibers of $\cF$.  In the case that
     $\cX$ is a tame, global quotient with projective coarse
     moduli space, we believe that there is again a \emph{global
     construction} of the flattening stratification $\Sigma\rightarrow
     S$ along with a partition labelled by the Hilbert polynomial (as
     defined in section~\ref{sec-HP}).  However we don't know what
     this global construction is, and the proof of the existence of
     $\Sigma\rightarrow S$ given above is the one truly
     non-constructive step in the proof of theorm~\ref{thm-thm2}.

\section{Hilbert polynomials}~\label{sec-HP}

Recall that a Deligne-Mumford stack $\cX$ is \emph{tame} if for
each algebraically-closed field $k$ and each 1-morphism $\zeta:\SP
k\rightarrow \cX$, the $k$-valued points of
\begin{equation}
G_{\zeta} := \SP k \times _{(\zeta,\zeta),\mc X\times \mc
X,\Delta} {\mc
  X}
\end{equation}
form a group of order  prime to char$(k)$.   We remind the reader
of some facts about tame Deligne-Mumford stacks.

\

\begin{lem}~\label{lem-tame}
Let $\cX$ be a tame Deligne-Mumford stack, $\pi:{\mc
  X}\rightarrow X$ its coarse moduli space.
\begin{enumerate}
\item The additive functor $\pi_*$ from the category of
  $\OO_{\cX}$-modules
  to the category of $\OO_X$-modules
  maps quasi-coherent sheaves to
  quasi-coherent sheaves and maps locally finitely-presented sheaves
  to locally finitely-presented sheaves.
\item The additive functor $\pi_*$ is exact, in particular
  $R^i\pi_*\cF=0$ for $i>0$ and $\cF$ any quasi-coherent $\OO_{\cX}$-module.
\item Suppose $g:X\rightarrow S$ is a morphism of algebraic spaces and
  suppose $\cF$ is a quasi-coherent sheaf on $\cX$ which is
  flat over $S$.  Then also $\pi_* \cF$ is flat over $S$.
\end{enumerate}
\end{lem}

\begin{proof}
  (1) and (2) form ~\cite[lemma 2.3.4]{AV}.  And (3) follows by the
  same local analysis in the proof of ~\cite[lemma 2.3.4]{AV} since
  the invariant submodule $M^\Gamma$ (in the notation of ~\cite{AV})
  of an $S$-flat module $M$ is a direct summand, and so it is also
  $S$-flat.
\end{proof}

Suppose $k$ is a field and $\cX$ is a separated, locally
finitely-presented tame Deligne-Mumford stack over $\SP{k}$, and
let $\pi :\cX\rightarrow X$ be its coarse moduli space. Suppose
that $\cF$ is a coherent $\OO_\cX$-module with proper support.
Then by lemma \ref{lem-tame}, $H^i_{\text{\'et}}(\cX, \cF) =
H^i_{\text{\'et}}(X, \pi _*\cF)$ for all $i$, and hence these
groups are zero for $i$ bigger than the dimension of $\cX$.
Therefore the sum
\begin{equation}
\chi(\cX,\cF) = \sum_{i=0}^\infty (-1)^i \text{dim}_k
H^i_{\text{\'et}}(\cX,\cF)
\end{equation}
is finite.  For each locally free sheaf $E$ of finite rank on
$\cX$, we also have that $E\otimes_{\OO_\cX} \cF$ is again
coherent with proper support. Since $E\mapsto \chi({\mc
X},E\otimes_{\OO_\cX} \cF)$ is additive in short exact sequences,
we have a well-defined group homomorphism
\begin{equation}
P_\cF:K^0(\cX)\rightarrow \ZZ, [E]\mapsto \chi({\mc
  X},E\otimes_{\OO_\cX} \cF).
\end{equation}

\

\begin{defn}~\label{defn-HP}
  Given a homomorphism of Abelian groups $i:A\rightarrow
  K^0(\cX)$, define the \emph{A-Hilbert polynomial} of $\cF$,
  $P_{A,\cF}$, to be the function $P_{A,\cF}=P_\cF\circ i$.
\end{defn}

\textbf{Remarks:} (1) If ${\mc L}$ is an invertible sheaf on $\cX$,
$A=\ZZ[x,x^{-1}]$ and $i:A\rightarrow K^0(\cX)$ is the group homomorphism
such that $i(x^d)=[{\mc L}^d]$, then the $A$-Hilbert polynomial of
$\cF$,
$P_{A,\cF}$ is the usual
Hilbert polynomial of $\cF$ with respect to ${\mc L}$.  We will
need to consider cases where $i:A\rightarrow K^0({\mc F})$ cannot be
reduced to this form, which is why our notion of Hilbert polynomial is
so general.

\

(2) The most instructive example, from our point of view, is when
    $\cX = BG$ for some finite, \'etale $k$-group scheme $G$.
    Then $K^0(BG)$ is naturally isomorphic to the Grothendieck group of the
    category of finite $k[G]$-modules, i.e. the representation ring of
    $k[G]$.  And the Hilbert polynomial $P_\cF$ determines the
    image $[\cF]$ of $\cF$ in $K^0(BG)$.

\

Now suppose that $f:\cX\rightarrow S$ is a 1-morphism from a
tame Deligne-Mumford stack to a connected algebraic space such that $f$ is
separated and locally finitely-presented,  Define $A=K^0(\cX)$ and
for each field $k$ and each morphism $g:\SP k\rightarrow S$, define
$i_g:A\rightarrow K^0(\SP k\times_S \cX)$ to be the pullback map
$K^0(\text{pr}_2)$.

\

\begin{lem}~\label{lem-HPconst}
Suppose that $\cF$ is a locally finitely-presented quasi-coherent
sheaf on $\cX$ which is flat over $S$ and which has proper support
over $S$.  Then there exists a function $P:A\rightarrow \ZZ$ such that
for all $g:\SP k\rightarrow S$, $P_{A,g^*\cF} = P$.
\end{lem}

\begin{proof}
We need to show that for each locally free sheaf $\cE$ on ${\mc
  X}$, the function
\begin{equation}
(g:\SP k\rightarrow S) \mapsto \chi\lt(g^*\lt(\cE\otimes_{\OO_{\mc
    X}}\cF\rt)\rt)
\end{equation}
is constant.  Since $R^i\pi_*$ vanishes on all quasi-coherent modules
for $i>0$, we have
\begin{equation}
\chi\lt(g^*\lt(\cE\otimes_{\OO_\cX}\cF\rt)\rt) =
\chi\lt(g^*\pi_*\lt( \cE\otimes_{\OO_\cX}\cF\rt)\rt).
\end{equation}
And by lemma~\ref{lem-tame} (3), we know that $\pi_*({\mc
  E}\otimes_{\OO_\cX}\cF)$ is an $S$-flat, locally
  finitely-presented sheaf with proper support over $S$.  So by
  ~\cite[corollary, p. 50]{Mum}, we conclude that there is some
  $P([\cE])\in \ZZ$ such that for all $g:\SP k\rightarrow S$, we
  have
\begin{equation}
\chi\lt(g^*\pi_*\lt( \cE\otimes_{\OO_\cX}\cF\rt)\rt) =
  P([\cE]).
\end{equation}
\end{proof}

Fix an additive homomorphism $P:K^0(\cX)\rightarrow \ZZ$ and
define $Q^P = Q^P(\cF/\cX/S)$ to be the subfunctor of $Q({\mc
  F}/\cX/S)$ such that for each $S$-scheme, $g:Z\rightarrow S$ we
have $Q^P(g:Z\rightarrow S)$ is the set of quotients $[g^*{\mc
  F}\rightarrow {\mc G}]\in Q(g:Z\rightarrow S]$ such that for each
field $k$ and each morphism $h:\SP k \rightarrow Z$, we have
$P_{A,h^*{\mc G}}=P$.  By lemma~\ref{lem-HPconst}, we see that $Q^P$ is
an open and closed subfunctor (possibly empty) of $Q$ and that $Q$ is
the disjoint union of $Q^P$ as $P$ ranges over all $P$.

\

Observe that theorem~\ref{thm-thm2} is implied by the following refinement.

\

\begin{thm}~\label{thm-thm3}
  Suppose that $S$ is an affine scheme.  Suppose that $\cX$
  is a tame Deligne-Mumford stack which is separated and
  finitely-presented over $S$, whose coarse moduli space is a
  quasi-projective $S$-scheme (resp. projective $S$-scheme), and which
  is a global quotient over $S$.
  Then for each locally
  finitely-presented quasi-coherent sheaf $\cF$ on $\cX$ and
  for each homomorphism $P:K^0(\cX)\rightarrow \ZZ$, the functor
  $Q^P(\cF/\cX/S)$ is represented by an algebraic space $Q^P$ which
  admits a factorization $Q^P\rightarrow Q'\rightarrow S$ where $Q'$
  is projective over $S$ and $Q^P\rightarrow Q'$ is a
  finitely-presented, quasi-finite monomorphism.
  If $\cF$ has proper support over $S$, then
  $Q^P\rightarrow Q'$ is a finitely-presented closed immersion.
\end{thm}

\textbf{Remark:}  In particular, if $S$ is affine, then
$Q^P(\cF/\cX/S)$ is represented by a quasi-projective $S$-scheme
(which is projective if the support of $\cF$ is proper over $S$).

\

We will prove theorem~\ref{thm-thm3} in sections \ref{gensheavsec} and
\ref{nattransec}.

\section{Generating sheaves}\label{gensheavsec}
Let $\cX$ be a tame Deligne-Mumford stack with coarse
moduli space
$\pi:\cX\rightarrow X$.  For each locally free sheaf $\cE$ on
$\cX$, define additive functors
\begin{eqnarray}
F_\cE:\text{Quasi-coherent}_\cX \rightarrow
\text{Quasi-coherent}_X, \\
G_\cE:\text{Quasi-coherent}_\cX \rightarrow
\text{Quasi-coherent}_\cX
\end{eqnarray}
by the formulas
\begin{equation}
F_\cE(\cF) = \pi_*\textit{Hom}_{\OO_{\cX}}(\cE,\cF), \ \ \ G_\cE =
\pi^*\lt(F_\cE(\cF)\rt)\otimes_{\OO_{\cX}}\cE
\end{equation}
where $\mc F$ is a quasi-coherent sheaf on $\mc X$.

\

By lemma~\ref{lem-tame}, $F_\cE$ is an exact functor which
preserves flatness and preserves the property of being locally
finitely-presented.  And $G_\cE$ is a right-exact functor which
preserves the property of being locally finitely-presented.  Moreover
there is a natural transformation $\theta_\cE:G_\cE\Rightarrow
\text{Id}$ where $\text{Id}$ is the identity functor on the category
of quasi-coherent sheaves on $\cX$ and for a quasi-coherent sheaf
$\cF$ on $\cX$, the morphism
\begin{equation}
\theta_{\cE}(\cF):\pi^* \lt( \pi_* \textit{Hom}_{\OO_{\cX}} (\cE,\cF)
\rt) \otimes_{\OO_{\cX}} \cE \rightarrow \cF
\end{equation}
is the left adjoint to the natural morphism
\begin{equation}
\pi^*\pi_*\textit{Hom}_{\OO_{\cX}}(\cE,\cF)\rightarrow
\textit{Hom}_{\OO_{\cX}}(\cE,\cF),
\end{equation}
which is itself the left adjoint of the identity morphism
\begin{equation}
\pi_* \textit{Hom}_{\OO_{\cX}} (\cE,\cF) \rightarrow \pi_*
\textit{Hom}_{\OO_{\cX}} (\cE,\cF).
\end{equation}

\textbf{Remark:} Since $\cX$ is a tame stack, $\pi_*$ from
quasi-coherent sheaves on $\cX$ to quasi-coherent sheaves on $X$
is compatible with \emph{arbitrary} base changes $T\rightarrow X$
(see for example the proof of ~\cite[2.3.3]{AV}).  It follows that
also the functors $F$, $G$ and the natural transformation $\theta$
are compatible with \emph{arbitrary} base changes $T\rightarrow
X$.

\begin{defn}~\label{defn-gen}  With notation as above, $\cE$ is a
  \emph{generator } for $\cF$ if $\theta_{\cE}(\cF)$ is surjective.
\end{defn}

\textbf{Example:}  Suppose that $G$ is a finite group and $\cX =
BG\times_{\SP \ZZ} \SP k$.  Then the quasi-coherent
$\OO_{\cX}$-modules correspond to (left-)modules over $k[G]$.  Let $\cE$ be
the locally free sheaf corresponding to the left regular
representation of $G$ on $k[G]$.  Then for every quasi-coherent sheaf
$\cF$, $\cE$ is a generator for $\cF$.  For that matter, if $M$ is any
$k[G]$-module which contains every irreducible representation of $G$
as a submodule, then the locally free sheaf corresponding to $M$ is a
generator for every $\cF$.

\

The goal of this section is to prove that the previous example is
typical for separated, tame, Deligne-Mumford stacks which are global
quotients.  As suggested by the referee, we first prove the following
\emph{pointwise} condition for a locally free sheaf $\cE$ to be a generator.

\begin{thm}~\label{thm-genpt}
Suppose that $\cX$ is a tame Deligne-Mumford stack as above and
$\cE$ is a locally free sheaf on $\cX$. For each
algebraically-closed field $k$ and each $1$-morphism $\zeta:\SP k
\rightarrow \cX$ with stabilizer group $G_k$ (c.f.
section~\ref{sec-HP}), let $\widetilde{\zeta}:BG_\zeta \rightarrow
\cX$ be the natural map.

Suppose that for some $1$-morphism $\zeta:\SP k \rightarrow \cX$,
the locally free sheaf $\cE_\zeta :=
\lt(\widetilde{\zeta}\rt)^*\cE$ is a generator for the left
regular representation of $G_\zeta$, i.e. considered as a
finite-dimensional $k[G_\zeta]$-module, $\cE_\zeta$ contains every
irreducible representation of $G_\zeta$ as a submodule. Then there
exists a Zariski open subset $U\subset X$ containing the image of
$\zeta(\SP k)$ such that the restriction of $\cE$ to $\cX \times_X
U$ is a generator for every quasi-coherent sheaf on $\cX\times_X
U$. In particular, if $\lt(\widetilde{\zeta} \rt)^*\cE$ is a
generator for the left regular representation of $G_\zeta$ for
every geometric point $\zeta:\SP k \rightarrow \cX$, then $\cE$ is
a generator for every quasi-coherent sheaf on $\cX$.
\end{thm}

\textbf{Remarks:} (1) The converse to theorem~\ref{thm-genpt} is
obvious.  Suppose that $\cE$ is a locally free sheaf on
$\cX$ and suppose that $\widetilde{\zeta}:BG_\zeta
\rightarrow \cX$ is as above.  Let $\cF$ be the push-forward by
$\widetilde{\zeta}$ of the sheaf on $BG_\zeta$ corresponding to the
left regular representation of $G_\zeta$.  Then $\cE$ is a generator
for $\cF$ iff $\lt(\widetilde{\zeta}\rt)^*\cE$ is a generator for the
left regular representation of $G_\zeta$.

\

(2) Suppose $\cX$ is a quotient stack of the form $\cX = [Y/GL_n]$.
    Then for each geometric point $\zeta:\SP k \rightarrow \cX$ the
    stabilizer group $G_\zeta$ is a closed subgroup scheme of
    $GL_{n,k}$ compatibly with the composite map
\begin{equation}
BG_\zeta \xrightarrow{\widetilde{\zeta}} \cX \rightarrow BGL_n.
\end{equation}
Of course every representation of $GL_n$ over $k$ induces a
representation of $G_\zeta$, i.e. every quasi-coherent sheaf on
$BGL_{n,k}$ pulls back to a quasi-coherent sheaf on $BG_\zeta$.
Moreover, it is well known that there exists a
$GL_{n,\ZZ}$-submodule
\begin{equation} W\subset
H^0(GL_{n,\ZZ},\OO_{GL_{n,\ZZ}})\end{equation} which is a
finitely-generated free $\ZZ$-module such that the composite map
\begin{equation}
W \otimes_\ZZ k \rightarrow H^0(GL_{n,k},\OO_{GL_{n,k}}) \rightarrow
H^0(G_\zeta, \OO_{G_\zeta})
\end{equation}
is surjective.  So for each point $\zeta:\SP k \rightarrow \cX$, there
is a $GL_{n,\ZZ}$ representation $W$ which is a finitely-generated
free $\ZZ$-module such that the pullback $\cE$ under $\cX
\rightarrow BGL_{n,\ZZ}$ of the locally free sheaf associated to $W$
satisfies the condition in theorem~\ref{thm-genpt}.  In order to apply
theorem~\ref{thm-genpt}, we need to find a single representation $W$
of $GL_{n,\ZZ}$ which works for every geometric point $\zeta:\SP k
\rightarrow \cX$.

\

(3) The stack $BG_{\zeta }$ is naturally isomorphic to the fiber
product
\begin{equation}
\cX \times _{X, \pi \circ \zeta }\SP k
\end{equation}
and the map $\widetilde {\zeta }$ is simply the projection onto
the first factor.

\

We will prove theorem~\ref{thm-genpt} after proving a proposition
about quotient stacks by a finite group.
Suppose that $\cE$ is a locally free sheaf on
$\cX$ which satisfies the hypotheses of theorem~\ref{thm-genpt} and
suppose that $\cF$ is a quasi-coherent sheaf on $\cX$.  To show that
$\cE$ is a generator for $\cF$, it suffices to prove this after
performing a faithfully-flat base change of the coarse moduli space
$X$.  By ~\cite[lemma 2.2.3]{AV}, there is an \'etale cover $(X_\alpha
\rightarrow X)$ such that for each $\alpha$ the fiber product $\cX
\times_X X_\alpha$ is of the form $[Y/G]$ where $Y$ is a finite scheme
over $X_\alpha$ and $G$ is a finite group which acts on $Y$ via
$X_\alpha$-morphisms.  Moreover, since $\cX$ is tame, we may assume
that the order of $G$ is prime to the characteristic of every residue
field of $X_\alpha$.  So the first step in proving
theorem~\ref{thm-genpt} is the case that $\cX = [Y/G]$ as above.

Now there is a map $f: [Y/G] \rightarrow BG$.  Let $E$ denote the
locally free sheaf on $BG$ corresponding to the left regular
representation of $G$.  If $\cE$ is a locally free sheaf on $\cX$
which is a generator for every quasi-coherent sheaf, then in
particular it is a generator for $f^*E$.  The following proposition
shows that this is a sufficient condition for $\cE$ to be a generator
for every quasi-coherent sheaf.

\begin{prop}~\label{prop-genBG}
Suppose that $\cX$ is a tame, separated Deligne-Mumford stack of the form
$\cX=[Y/G],$ as above.  Let $f:\cX \rightarrow BG$ and $E$ be as above.
Then $f^*E$ is a
generator for every quasi-coherent $\OO_{\cX}$-module $\cF$.
\end{prop}

\begin{proof}
Let $g:Y\rightarrow \cX$ denote the quotient $1$-morphism, and let
$p:\cX\rightarrow X$ denote the map to the coarse moduli space of
$\cX$.
Observe
that $f^* E$ is simply $g_*\OO_Y$.
The
composition $p\circ g:Y\rightarrow X$ is a finite, surjective morphism
of algebraic spaces.  In particular, it is affine and for each
quasi-coherent $\OO_Y$-module $\cG$, the induced morphism
\begin{equation}
\alpha:(p\circ g)^*(p\circ g)_* \cG \rightarrow \cG,
\end{equation}
is surjective.  By adjointness of $p_*$ and $p^*$, there is an induced
morphism
\begin{equation}
p^* (p\circ g)_* \cG \rightarrow g_* \cG,
\end{equation}
in fact this is precisely $\theta_{\OO_{\cX}}(g_*\cG)$.
Since $g_*\cG$ is a module over $g_*\OO_Y$, we get an induced
morphism of $\OO_{\cX}$-modules
\begin{equation}
\phi:g_*\OO_Y \otimes_{\OO_{\cX}} p^*(p\circ g)_* \cG \rightarrow g_*\cG.
\end{equation}
The claim is that $\phi$ is surjective; let us assume this for a
moment.  The canonical injection $g^\#:\OO_{\cX}\rightarrow
g_*\OO_Y$ induces a morphism of $\OO_{\cX}$-modules
\begin{equation}
\psi:g_*\OO_Y \otimes_{\OO_{\cX}}p^*
p_*\textit{Hom}_{\OO_{\cX}}(g_*\OO_Y,g_*\cG) \rightarrow
g_*\OO_Y \otimes_{\OO_{\cX}}p^*
p_*\textit{Hom}_{\OO_{\cX}}(\OO_{\cX},g_*\cG).
\end{equation}
Since $g$ is representable, finite and flat, there is a surjective
trace map
\begin{equation}
t:g_*\OO_Y\rightarrow \OO_{\cX},
\end{equation}
which splits the injection $g^\#$.
Therefore $\psi$ is surjective.  And
$\theta_{\cE}(g_*\cG)$ is the composite $\phi\circ\psi$.  Since $\psi$
and $\phi$ are both surjective, we conclude that
$\theta_{\cE}(g_*\cG)$ is surjective.  So to prove that
$\theta_{\cE}(g_*\cG)$ is surjective, it suffice to prove the claim
that $\phi$ is surjective.

\

To see that $\phi$ is surjective, we apply $g^*$; since $g$ is
faithfully flat, we may check surjectivity after base-change by $g$.
Since $g$ is affine, the canonical morphism
\begin{equation}
 (g^* g_* \OO_Y)\otimes_{\OO_Y}\cG\rightarrow g^*g_*\cG
\end{equation}
is an isomorphism.
Similarly, the canonical morphism
\begin{equation}
(g^* g_* \OO_Y)\otimes_{\OO_Y} (p\circ g)^*(p\circ g)_*\cG\rightarrow
g^*\lt(g_*\OO_Y \otimes_{\OO_{\cX}} p^*(p\circ g)_* \cG\rt)
\end{equation}
is also an isomorphism.
And we have a commutative diagram
\begin{equation}
\begin{CD}
(g^* g_* \OO_Y)\otimes_{\OO_Y} (p\circ g)^*(p\circ g)_*\cG @>
1\otimes\alpha >> (g^* g_* \OO_Y)\otimes_{\OO_Y}\cG \\
@VVV  @VVV \\
g^*\lt(g_*\OO_Y \otimes_{\OO_{\cX}} p^*(p\circ g)_* \cG\rt) @>g^*\phi
>> g^*g_*\cG
\end{CD}
\end{equation}
Since $\alpha$ is surjective, so is $1\otimes \alpha$.  Therefore
$g^*\phi$ is surjective, and it follows that $\phi$ is surjective.
This proves the claim, and we conclude that $f^* E$ is a generator for
all sheaves of the form $g_*\cG$ with $\cG$ a quasi-coherent
$\OO_Y$-module.

\

Let $\cF$ be a quasi-coherent $\OO_{\cX}$-module.
Since $g$ is affine, the
canonical morphism of $\OO_{\cX}$-modules
\begin{equation}
g_*\OO_Y \otimes_{\OO_{\cX}} \cF \rightarrow g_* g^* \cF,
\end{equation}
is an isomorphism.  Therefore
$f^* E$ is a generator for $g_*\OO_Y \otimes_{\OO_{\cX}} \cF$.
And we have a surjective morphism of $\OO_{\cX}$-modules,
\begin{equation}
t\otimes 1: g_*\OO_Y \otimes_{\OO_{\cX}} \cF \rightarrow \cF.
\end{equation}
Since $G_{f^*E}$ is a right-exact functor, we conclude that $f^* E$ is
also a generator for $\cF$.
\end{proof}

\begin{cor}~\label{cor-genBG}
Suppose that $\cX$, $f:\cX \rightarrow BG$, and $E$ are as in
proposition~\ref{prop-genBG}.  If $\cE$
is a locally free sheaf on $\cX$ which generates $f^* E$, then $\cE$
generates $\cF$ for every quasi-coherent $\OO_{\cX}$-module $\cF$.
\end{cor}

\begin{proof}
Since $p_*$ is exact, for any quasi-coherent $\OO_X$-module $\cG$ and
any quasi-coherent $\OO_{\cX}$-module $\cF$, we have that
\begin{equation}
p_*(p^*\cG\otimes_{\OO_{\cX}} \cF) = \cG \otimes_{\OO_{X}} p_*\cF.
\end{equation}
To see this, it suffices to work locally over affine opens in $X$, so
we may suppose that $\cG$ is a colimit of finitely-presented
$\OO_X$-modules.  Since $p_*$ commutes with colimits, we are reduced
to the case that $\cG$ is finitely-presented.  Since $p^*$ is
right-exact and since $p_*$ is exact, we are reduced to the case that
$\cG=\OO_X$, which is trivial.

\

Similarly, we conclude that $\theta_{\cE}(p^*\cG \otimes_{\OO_{\cX}}
\cF)$ is just $1_{p^* \cG}\otimes \theta_{\cE}(\cF)$.  In particular,
if $\cE$ generates $\cF$, then $\cE$ generates
$p^*\cG\otimes_{\OO_{\cX}} \cF$.  Therefore,  $\cE$ generates every
sheaf of the
form $p^*\cG\otimes_{\OO_{\cX}} f^* E$.  But by
proposition~\ref{prop-genBG}, we conclude that every quasi-coherent
sheaf $\cF$ is a surjective image of a sheaf of the form
$p^*\cG\otimes_{\OO_{\cX}} f^*E$.  Since $G_{\cE}$ is right-exact, we
conclude that $\cE$ generates $\cF$.
\end{proof}

Now we come to the proof of theorem~\ref{thm-genpt}:

\begin{proof}
Suppose $\cE$ satisfies the hypotheses of theorem~\ref{thm-genpt}.
By ~\cite[lemma 2.2.3]{AV}, there is an \'etale morphism $g:X'
\rightarrow X$ whose image contains the image of $\zeta(\SP k)$ and
such that $\cX' = \cX \times_X X'$ is of the form $[Y/G]$ for some
finite $X'$-scheme $Y$, for some finite group $G$ and an action of $G$
on $Y$ via $X'$-morphisms.  Let $f:[Y/G] \rightarrow BG$ and $E$ be as
in proposition~\ref{prop-genBG}.  Let $\cE'$ denote the pullback of
$\cE$ to $\cX'$.  Consider the morphism $\theta_{\cE'}(f^* E)$.
Denote the cokernel of this morphism by $\cG$.  Then $\cG$ is a
finitely-generated quasi-coherent sheaf, and therefore has closed
support $Z\subset \cX$.

Since $k$ is algebraically-closed and since
$X'\rightarrow X$ is finite, we can find a factorization $\xi:\SP k
\rightarrow X'$ of $\SP k \rightarrow X$.
Now the functors $F$ and $G$
and the natural transformation $\theta$ are compatible with arbitrary
base change on $X$.  Therefore the pullback of $\theta_{\cE'}(f^* E)$
equals $\theta_{\xi^*\cE'}(\xi^* f^* E)$.  Of course we have
\begin{equation}
\cX' \times_{X'} \SP k = \cX \times_X \SP k = BG_\zeta.
\end{equation}
By hypothesis $\xi^* \cE' = \widetilde{\zeta}(\cE)$ generates the
regular representation of $G_\zeta$.  So by
corollary~\ref{cor-genBG}, we conclude that
$\theta_{\xi^*\cE'}(\xi^* f^* E)$ is surjective.  So the image of
$\xi(\SP k)$ is not in the support of $\cG$.  Denoting by
$V\subset X'$ the complement of the image of the support of $\cG$,
we have that $V$ contains $\xi(\SP k)$ and the restriction of
$\cE'$ over $V$ generates $f^* E$.

Let $U$ be the image of $V\rightarrow X$ so that the canonical
morphism $h:\cX \times_X V
\rightarrow \cX \times_X U$ is faithfully flat and $U$ contains the
image of $\zeta(\SP k)$.
Let $\cF$ be any quasi-coherent sheaf on $\cX \times_X U$.  Then
$h^*\theta_{\cE}(\cF)$ is $\theta_{h^*\cE}(h^*\cF)$.  By
corollary~\ref{cor-genBG}, we know that $\theta_{h^*\cE}(h^*\cF)$ is
surjective.  Since $h$ is faithfully flat, we conclude that also
$\theta_{\cE}(\cF)$ is surjective.  So the restriction of $\cE$ to
$\cX \times_X U$ is a generator for every quasi-coherent sheaf $\cF$
on $\cX \times_X U$.  This completes the proof of the theorem.
\end{proof}

Now we come to the main theorem of this section:

\begin{thm}~\label{thm-genthm}
Suppose that $\cX$ is a quasi-compact, tame Deligne-Mumford stack of
the form $[Y/GL_{n,\ZZ}]$ for some scheme $Y$ and an action of
$GL_{n,\ZZ}$.  Let $f:\cX \rightarrow BGL_{n,\ZZ}$ be the obvious
$1$-morphism.
There exists a finitely-generated free
$\ZZ$-module, $W$ which is a representation of $GL_{n,\ZZ}$ such that
if $E$ on $BGL_{n,\ZZ}$ is the induced locally free sheaf then $f^* E$
has the following property:  For every morphism of algebraic spaces
$g:X'\rightarrow X$ with corresponding map $h:\cX \times_X X'
\rightarrow \cX$ and for every quasi-coherent sheaf $\cF$ on $\cX
\times_X X'$, we have that $h^*f^* E$ is a generator for $\cF$.
\end{thm}

\begin{proof}
First of all suppose we prove the result when $g:X'\rightarrow X$ is
the identity map.  Then by remark $(1)$ following the statement of
theorem~\ref{thm-genpt}, it follows that for every
algebraically-closed field $k$ and every $1$-morphism $\zeta:\SP k
\rightarrow \cX$
we have that $\widetilde{\zeta}^*f^* E$ generates the left regular
representation of $G_\zeta$.
Now suppose that $g:X'\rightarrow X$ is an arbitrary morphism, $k$ is
an algebraically-closed field and $\xi:\SP k \rightarrow \cX \times_X
X'$ is a $1$-morphism.  Let $\zeta:\SP k \rightarrow \cX$ be $h\circ
\xi$.  Then $G_\xi = G_\zeta$, and $\xi^*h^*f^*E = \zeta^* f^* E$.  So
we see that $h^* f^*E$ on $\cX \times_X X'$ satisfies the hypotheses
of theorem~\ref{thm-genpt}.  So $h^* f^* E$ generates every
quasi-coherent sheaf.  Therefore to prove our theorem, we need only
consider the case that $g:X'\rightarrow X$ is the identity map.

By remark $(2)$ following the statement of theorem~\ref{thm-genpt},
for each geometric point $\zeta:\SP k \rightarrow \cX$, there is a
finitely-generated free $\ZZ$-module $W_\zeta$ which is a
representation of $GL_{n,\ZZ}$ and such that if $E_\zeta$ is the
corresponding locally free sheaf on $BGL_{n,\ZZ}$, then
$\widetilde{\zeta}^* f^* W_\zeta$ generates the left regular
representation.  Then by theorem~\ref{thm-genpt}, there is a Zariski
open neighborhood $U_\zeta \subset X$ of the image of $\zeta(\SP k)$
such that the restriction of $f^* E_\zeta$ to $\cX \times_X U_\zeta$
generates all quasi-coherent sheaves.  Allowing $k$ and $\zeta$ to
vary over all algebraically-closed fields and $1$-morphisms, we obtain
a Zariski open cover $\{ U_\zeta\}$ of $X$.  Since $\cX$ is
quasi-compact, some finite subset covers $X$, say
$\zeta_1,\dots,\zeta_N$.  Let $W$ be the $GL_{n,\ZZ}$-representation
$W_{\zeta_1} \oplus \dots \oplus W_{\zeta_N}$.  Then clearly $W$
satisfies the conditions in the theorem.
\end{proof}

\section{Natural transformation of Quot functors}\label{nattransec}
Suppose that $S$ is a quasi-compact algebraic space, $f:\cX\rightarrow
S$ is a separated, finitely-presented $1$-morphism of a tame
Deligne-Mumford stack to $S$ such that $\cX$ is a global quotient over
$S$, say $\cX=[Y/GL_{n,X}]$.  By theorem~\ref{thm-genthm}, there is a
$GL_{n,\ZZ}$ representation $V$ which is a finite free $\ZZ$-module
and such that for every morphism of algebraic spaces $T\rightarrow S$
and every quasi-coherent sheaf $\cG$ on $T\times_S \cX$, $\cE_V$
generates $\cG$.

\

Suppose that $\cF$ is a locally finitely-presented $\OO_{\cX}$-module.
Suppose that $P:K^0(\cX)\rightarrow \ZZ$ is a given Hilbert
polynomial.
Let $p:\cX\rightarrow X$ denote the coarse moduli scheme of $\cX$ and
let $P_V:K^0(X)\rightarrow \ZZ$ denote the map
\begin{equation}
P_V([E]) = P([\textit{Hom}_{\OO_{\cX}}(\cE_V,p^*E)]).
\end{equation}
Suppose $T\rightarrow S$ is a morphism of a scheme to $S$ and suppose
that $\cF_T\rightarrow \cG$ is an element of $Q^P(\cF/\cX/S)(T)$.
Consider $F_{\cE_V}(\cF_T)\rightarrow F_{\cE_V}(\cG)$.  Since
$F_{\cE_V}$ is exact, this is still surjective.  By
lemma~\ref{lem-tame}, $F_{\cE_V}(\cG)$ is locally finitely-presented,
has proper support over $T$ and is flat over $T$.  Moreover every
geometric fiber has Hilbert polynomial $P_V$,
i.e. $F_{\cE_V}(\cF_T)\rightarrow F_{\cE_V}(\cG)$ is an element of
$Q^{P_V}(F_{\cE_V}/X/S)(T)$.  This defines a natural transformation
\begin{equation}
F_{\cE_V}:Q^P(\cF/\cX/S) \Rightarrow Q^{P_V}(F_{\cE_V}(\cF)/X/S).
\end{equation}

\begin{lem}~\label{lem-mono}
The natural transformation $F_{\cE_V}$ is a monomorphism, i.e. for
each $T\rightarrow S$ the morphism of sets
\begin{equation}
F_{\cE_V}(T):Q^P(\cF/\cX/S)(T)\rightarrow
Q^{P_V}(F_{\cE_V}(\cF)/X/S)(T),
\end{equation}
is an injection of sets.
\end{lem}

\begin{proof}
Given any morphism $T\rightarrow S$ of a scheme to $S$ and given any
element $\alpha:F_{\cE_V}(\cF)_T\rightarrow G$ in
$Q^{P_V}(F_{\cE_V}(\cF)/X/S)(T)$, define $\beta:K\rightarrow
F_{\cE_V}(\cF)_T$ to be the kernel of
$\alpha$ and define $\eta_T(\alpha):\cF\rightarrow \cG$ to be the
cokernel of the composition:
\begin{equation}
\begin{CD}
p^*K\otimes_{\OO_{\cX}} \cE_V @> p^*\beta\otimes 1 >> G_{\cE_V}(\cF_T)
@> \theta_{\cE_V}(\cF_T) >> \cF_T.
\end{CD}
\end{equation}

\

Now suppose given $\nu:\cF_T\rightarrow \cG$ in $Q^P(\cF/\cX/S)(T)$
and let $\mu:\cK\rightarrow \cF_T$ denote the kernel of $\nu$.  Then
we have a short exact sequence:
\begin{equation}
\begin{CD}
0 @>>> F_{\cE_V}(\cK) @> F_{\cE_V}(\mu) >> F_{\cE_V}(\cF)_T @>
F_{\cE_V}(\nu) >> F_{\cE_V}(\cG) @>>> 0.
\end{CD}
\end{equation}
Notice that if $\alpha=F_{\cE_V}(\nu)$, then $\beta=F_{\cE_F}(\mu)$.
Since $p^*$ is right-exact, we have a commutative diagram of exact
sequences:
\begin{equation}
\begin{CD}
& & G_{\cE_V}(\cK) @> p^*\beta\otimes 1 >>
    G_{\cE_V}(\cF_T) @> p^*\alpha\otimes 1 >>
      G_{\cE_V}(\cG) @>>> 0 \\
& & @V \theta_{\cE_V}(\cK) VV     @V \theta_{\cE_V}(\cF_T) VV
     @V \theta_{\cE_V}(\cG) VV \\
0 @>>> \cK @> \mu >> \cE_T @> \nu >> \cK @>>> 0
\end{CD}
\end{equation}
It follows by the snake lemma that we have the formula:
\begin{equation}~\label{eqn-main}
\eta_T(F_{\cE_V}(\mu)) = \mu.
\end{equation}
Therefore $F_{\cE_V}$ is injective.
\end{proof}

\textbf{Remark:}  Notice that the association $T\mapsto \eta_T$ is
functorial for arbitrary $S$-morphisms $T_1\rightarrow T_2$.  This
might seem odd since the formation of the kernel $K$ is not compatible
arbitrary base-change.  But one could define $\eta_T$ equivalently to
be the coequalizer of the diagram:
\begin{equation}
\begin{CD}
G_{\cE_V}(\cF_T) @> p^*\alpha \otimes 1 >> p^* G \otimes_{\OO_{\cX}}
\cE_V \\
@V \theta_{\cE_V}(\cF_T) VV \\
\cF_T
\end{CD}
\end{equation}
Since coequalizers are compatible with arbitrary base-change, we see
that $\eta$ is compatible with arbitrary $S$-morphisms $T_1\rightarrow
T_2$.

\

\begin{prop}~\label{prop-relrep}
The monomorphism $F_{\cE_V}$ of functors is relatively
representable by schemes.  In fact, $F_{\cE_V}$ is a
finitely-presented, quasi-finite, monomorphism.  If $\cF$ has
proper support over $S$, then $F_{\cE_V}$ is a finitely-presented,
finite, monomorphism, i.e. a finitely-presented closed immersion.
\end{prop}

\begin{proof}
Suppose $T\rightarrow S$ is a morphism of a scheme to $S$ and suppose
that $\alpha:F_{\cE_V}(\cF_T)\rightarrow G$ is an element of
$Q^{P_V}(F_{\cE_V}(\cF)/X/S)$.  Now form
$\eta_T(\alpha):\cF_T\rightarrow \cG$.  By theorem~\ref{thm-fs}, there
is a flattening stratification $g:\Sigma\rightarrow T$ for $\cG$ and
$g$ is a surjective, finitely-presented, quasi-finite monomorphism.
By lemma~\ref{lem-HPconst}, for each connected component $\Sigma_i$ of
$\Sigma$, the restriction of $\cG$ to $\Sigma_i$ has constant Hilbert
polynomial $P_i$.  In particular, there is a connected component
(possibly empty) $\Sigma_i$ of $\Sigma$ on which the restriction of
$\cG$ has Hilbert polynomial $P$.  Of course we have an induced
natural transformation of functors
\begin{equation}
\Sigma_i\rightarrow T\times_{\alpha,Q^{P_V}(F_{\cE_V}(\cF)/X/S)}
Q^P(\cF/\cX/S).
\end{equation}
By equation~\ref{eqn-main}, we have an inverse natural
transformation.  Thus we conclude that the fiber functor
$T\times_{\alpha,Q^{P_V}(F_{\cE_V}(\cF)/X/S)} Q^P(\cF/\cX/S)$ is
represented by the morphism $g_i:\Sigma_i\rightarrow T$.
Notice that this is a finitely-presented, quasi-finite monomorphism of
schemes.

\

In case $\cF$ has proper support over $S$, we know from
theorem~\ref{thm-thm1} that $Q^P\rightarrow S$ satisfies the valuative
criterion of properness.  Therefore $g_i:\Sigma_i\rightarrow T$
satisfies the valuative criterion of properness, i.e. $g_i$ is
finite.  But a finite monomorphism is precisely a closed immersion,
therefore $g_i$ is a finitely-presented closed immersion.
\end{proof}

Now we are in a position to prove theorem~\ref{thm-thm3}.  By
proposition~\ref{prop-relrep}, we know that
\begin{equation}
Q^P(\cF/\cX/S) \rightarrow Q^{P_V}(F_{\cE_V}(\cF)/X/S)
\end{equation}
is relatively representable by a
finitely-presented, quasi-finite monomorphism
(resp. finitely-presented closed immersion).

\

Now the proof that $Q^{P_V}(F_{\cE_V}(\cF)/X/S)\rightarrow S$ is
represented by a scheme which is quasi-projective over $S$
is essetially ~\cite[th\'eor\`eme 3.2, part IV]{FGA}.  In his proof,
Grothendieck makes some Noetherian hypotheses which are eliminated in
a standard way, cf.~\cite[section IV.8.9]{EGA}.  For the sake of
completeness, we include the proof here.

\

We have a locally closed immersion $X\hookrightarrow \PP^N_S$ for some
$N$.
By ~\cite[proposition IV.8.9.1]{EGA}, we can find a finite-type affine
$\ZZ$-scheme, $S'$, a quasi-projective (resp. projective) $S'$-scheme,
$X'$, and a coherent sheaf $F'$ on $X'$ along with a morphism
$S\rightarrow S'$ such that $X$ is isomorphic to $S\times_{S'} X'$,
and under this isomorphism $F_{\cE_V}(\cF)$ is isomorphic to the
pullback of $F'$.  Let $p:\ZZ[t,t^{-1}]\rightarrow \ZZ$ be the
polynomial $p(t^n) = P_V([\OO_X(n)])$ where $\OO_X(n)$ is the pullback
to $X$ of the invertible sheaf $\OO_{\PP^N_S}(n)$ on $\PP^N_S$.
polynomial, and consider the Quot functor $Q^p(F'/X'/S')$
which is the connected component of $Q(F'/X'/S')$ parametrizing
families of quotients $F'_T\rightarrow G$ such that for each closed
point $t\in T$, we have $\chi(G_t\otimes \OO(n)) = p(t^n)$.  By
lemma~\ref{lem-HPconst}, it follows
that $Q^{P_V}(F_{\cE_V}(\cF)/X/S)$ is isomorphic to a connected
component of the fiber product $S\times_{S'} Q^p(F'/X'/S')$.  So to
prove that $Q^{P_V}(F_{\cE_V}(\cF)/X/S)\rightarrow S$ is
quasi-projective, it suffices
to show that $Q^p(F'/X'/S')\rightarrow S'$ is quasi-projective.

\

By ~\cite[th\'eor\`eme 3.2, part IV]{FGA}, we know that the
restriction of the functor $Q^p(F'/X'/S')$ to the category of
locally Noetherian $S'$-schemes is represented by a
quasi-projective $S'$-scheme.  Using ~\cite[proposition
IV.8.9.1]{EGA} and ~\cite[th\'eor\`eme IV.11.2.6]{EGA}, it follows
that in fact this quasi-projective $S'$-scheme represents
$Q^p(F'/X'/S')$ on the category of all $S'$-schemes.  Thus
$Q^{P_V}(F_{\cE_V}(\cF)/X/S)\rightarrow S$ is represented by a
quasi-projective scheme.  This completes the proof of
theorem~\ref{thm-thm3}.

\bibliography{my}
\bibliographystyle{abbrv}

\end{document}